\documentclass[smallextended]{svjour3}     % onecolumn (second format)
\smartqed  % flush right qed marks, e.g. at end of proof
\usepackage{graphicx}
\usepackage{amssymb,amsfonts,amsmath,latexsym,mathrsfs}

\newcommand{\be}{\begin{enumerate}}
\newcommand{\ee}{\end{enumerate}}
\newcommand{\bm}[1]{{\mbox{\boldmath $#1$}}}
\newcommand{\beq}{\begin{equation}}
\newcommand{\eeq}{\end{equation}}
\newcommand{\beas}{\begin{eqnarray*}}
\newcommand{\eeas}{\end{eqnarray*}}

\newcommand{\zz}{\mathbb{Z}}
\newcommand{\rr}{\mathbb{R}}
\newcommand{\pp}{\mathbb{P}}
\newcommand{\hnq}{{\cal H}_n(q)}
\newcommand{\sn}{\mathfrak{S}_n}
\newcommand{\inv}{\mathrm{inv}}

\begin{document}

\title{Some Hecke Algebra Products and Corresponding Random
  Walks\thanks{Rosena Du is partially supported by the
    National Science Foundation of China under Grant No.~10801053 and
    No.~10671074. Research carried out when she was a Visiting Scholar
    at M.I.T.\ during the 2007--2008 academic year. Richard Stanley's
   contribution is based upon work supported by the National
 Science Foundation under Grant No.~DMS-0604423. Any opinions,
 findings and conclusions or recommendations expressed in this
 material are those of the author and do not necessarily reflect those
 of the National Science Foundation.}}
%General acknowledgments should be placed at the end of the article.{}


\titlerunning{Hecke Algebra Products}        % if too long for running head

\author{Rosena R. X. Du        \and
        Richard P. Stanley %etc.
}


\institute{           Rosena R. X. Du \at
              Department of Mathematics, East China
    Normal University, Shanghai 200041 \\
     \email{rxdu@math.ecnu.edu.cn}
            \and
              Richard P. Stanley \at
              Department of Mathematics, M.I.T., Cambridge, MA 02139 \\
              Tel.: 617-253-7930\\
              Fax: 617-253-4358\\
              \email{rstan@math.mit.edu}           %  \\
}

\date{Received: date / Accepted: date}

\maketitle

\begin{abstract}
 Let $\bm{i}=1+q+\cdots+q^{i-1}$. For certain
sequences $(r_1,\dots,r_l)$ of positive integers, we show that in
the Hecke algebra $\hnq$ of the symmetric group $\sn$, the product
$(1+\bm{r_1}T_{r_1})\cdots (1+\bm{r_l}T_{r_l})$ has a simple
explicit expansion in terms of the standard basis $\{T_w\}$. An
interpretation is given in terms of random walks on $\sn$.
\end{abstract}

\section{Main Results}

Let $\mathfrak{S}_n$ denote the symmetric group of permutations of
$\{1,\ldots, n\}$. For $1\leq i\leq n-1$ let $s_i=(i,i+1)\in\sn$,
the \emph{adjacent transposition} interchanging $i$ and $i+1$ and
leaving all the other elements fixed. For any $w \in \mathfrak{S}_n$
denote by $\ell(w)$ the \emph{length} of $w$, i.e., the minimal $p$
such that $w$ can be written as
  $$ w=s_{r_1}s_{r_2}\cdots s_{r_p} $$
for certain $r_1, r_2, \ldots, r_p$; such a sequence $r=(r_1,\ldots,
r_p)$ is called a \emph{reduced decomposition} (or reduced word)
provided $p=\ell(w)$.

The \emph{Hecke Algebra} (or \emph{Iwahori-Hecke algebra}) $\hnq$ of
the symmetric group $\sn$ (e.g., \cite[{\S}7.4]{JHumphreys}) is
defined as follows: $\hnq$ is an $\rr$-algebra with identity $1$ and
generators $T_1, T_2, \ldots, T_{n-1}$ which satisfy relations
\begin{eqnarray}
(T_i+1)(T_i-q)&=&0, \nonumber\\
T_i T_j&=&T_j T_i, \ \ \ |i-j|\geq 2,\\
T_i T_{i+1} T_i&=&T_{i+1} T_{i} T_{i+1}, \ \ \ 1\leq i \leq n-2.
\nonumber
\end{eqnarray}

For any $w =s_{r_1}s_{r_2}\cdots s_{r_p} \in \sn$ for which $(r_1,
r_2, \dots, r_p)$ is reduced, define $T_w =T_{r_1}T_{r_2}\cdots
T_{r_p}$. A basic property of Hecke algebras is that $T_w$ does not
depend on the choice of reduced decomposition of $w$, and for $1\leq
k \leq n-1$, $T_w$ satisfies

\begin{equation}\label{Tprod}
T_wT_{k}=\left\{
\begin{array}{ll}
T_{w s_k}, & \mbox{if}\  \ell(w s_k)=\ell(w)+1, \\
q T_{w s_k}+(q-1)T_w, & \mbox{if}\  \ell(w s_k)=\ell(w)-1.
\end{array}
\right.
\end{equation}

Let $r=(r_1, r_2, \cdots, r_l)$ be any sequence of positive integers
(not necessarily reduced).  For convenience assume that
$\max\{r_1,\dots,r_l\}=n-1$. Set
  $$ \bm{i}=1+q+\cdots+q^{i-1}. $$
For any $w \in \sn$, define $\alpha_r(w)\in\zz[q]$ by
\begin{equation*}
Q(r):=(1+\bm{r_1}T_{r_1})(1+\bm{r_2}T_{r_2})\cdots
(1+\bm{r_l}T_{r_l})=\sum_{w \in\sn} \alpha_r(w)T_w.
\end{equation*}
We are primarily concerned with the polynomials $\alpha_r(w)$. In
particular, for which $r$'s will $\alpha_r(w)$ have ``nice'' values
for all $w\in\sn$?

For each $w=w_1 w_2 \cdots w_n \in \sn$, we write $w \preceq r$ if
$w=s_{c_1}\cdots s_{c_k}$ for some subsequence $c_1,\dots,c_k$ of
$r=(r_1,\dots,r_l)$. This defines the \emph{Bruhat order} on $\sn$
\cite[{\S}5.9]{JHumphreys}. It follows from equation~\eqref{Tprod} that
$\alpha_r(w)=0$ unless $w\preceq r$.  Let $a_r(i)$ denote the number
of $i$'s in $r$, and let $\inv(w)=(\inv_w(1),\inv_w(2), \ldots,
\inv_w(n-1))$ denote the \emph{inversion sequence} of $w$, i.e., for
any $1 \leq i \leq n-1$, $\inv_w(i)$ is the number of $j$'s such
that $w_j<w_i$ and $j>i$.

We recursively define a sequence $r=(r_1,\dots,r_l)$ to be a
\emph{tight sequence} if it satisfies the following:
\be
\item $r_1=1$;

\item If $r$ is a tight sequence, and $k$ is a positive integer,
then $r'=(r,k)$ (the concatenation of $r$ and $k$) is also a tight
sequence if $ a_r(k)\leq a_r(k-1)-1$, and equality holds when there
exists a permutation $w$ satisfying $w \preceq r'$ but $w \npreceq
r$. \ee

For example, there are $6$ tight sequences of length $4$: 1111,
1211, 1231, 1212, 1213, 1234. And any prefix of the sequences
$(1,2,1,3,2,1,\ldots)$ or $(1,2,\ldots, n,$
$1,2,\ldots,n-1,\ldots,1,2,1)$ is a tight sequence.

The main result of this paper is the following.

\begin{theorem}\label{mainthm}
Let $r$ be a tight sequence with $\max\{r\}=n-1$. Then for any $w
\in \sn$ and $w \preceq r$, we have
\begin{equation} \label{main0}
\alpha_r(w)=\prod_{i=2}^{n-1} \bm{i}^{\max\{a_r(i-1)-1,
\inv_w(i)\}}.
\end{equation}
\end{theorem}

\begin{example}
\be\item[(a)] Define the \emph{standard tight sequence} $\rho_n$ of
degree $n$ by
   $$ \rho_n = (1,2,1,3,2,1,4,3,2,1,\dots,n-1,n-2,\dots,1). $$
It is easy to see that $\rho_n$ is not only a tight sequence but
also a reduced decomposition of the element
$w_0=n,n-1\dots,1\in\sn$. Theorem~\ref{mainthm} becomes
   \beq \alpha_{\rho_n}(w) = \bm{1}^{n-1}\bm{2}^{n-2}\cdots
           (\bm{n-1}), \label{eq:qunif} \eeq
independent of $w\in\sn$.
  \item[(b)] Let $r=(1,2,1,1,3,1)$. Then we have
$$\alpha_{r}(w)=\bm{2}^{3}, \ \ \forall
w\in\{1234,1324,2134,2314,3124,3214\}$$ and
$$\alpha_{r}(w)=\bm{2}^{3}\bm{3},\ \ \forall
w\in\{1243,1342,2143,2341,3142,3241\}.$$ Otherwise we have
$\alpha_r(w)=0$. \ee
\end{example}

Alexander Molev has pointed out (private communication dated
September 1, 2008) that Theorem~\ref{mainthm} in the case of the
standard tight sequence is connected to the ``fusion procedure'' for
the Hecke algebra, which goes back to Cherednik
\cite{cher1}\cite{cher2}.

\section{Proof of the Main Theorem}

For the proof of Theorem~\ref{mainthm} we need the following lemma.

\begin{lemma}\label{lem0}
Let $r=(r_1,\dots,r_l)\in\pp^l$ with $\max\{r\}=n-1$, where
$\pp=\{1,2,\dots\}$. Set $r'=(r,k)$, $1\leq k \leq n-1$. Then for
any $w \in \sn$ and $w \preceq r'$, we have \be
\item If $w\npreceq r$, then
$\alpha_{r'}(w)=\alpha_{r}(ws_k)\cdot \bm{k}$,

\item If $w \preceq r, w s_k \npreceq r$, then
$\alpha_{r'}(w)=\alpha_{r}(w)$,

\item If $w,w s_k \preceq r$, and $\ell(w s_k)=\ell(w)+1$, then
$\alpha_{r'}(w)= \alpha_{r}(w)+ \alpha_{r}(ws_k)\cdot \bm{k}q$,

\item If $w,w s_k \preceq r$, and $\ell(w s_k)=\ell(w)-1$, then
$\alpha_{r'}(w)=\alpha_{r}(w)\cdot q^k+ \alpha_{r}(ws_k)\cdot
\bm{k}$. \ee
\end{lemma}

\proof We have
\begin{equation}\label{indLemma}
Q(r')=Q(r)(1+\bm{k}T_k)
%=\sum_{\sigma \preceq r}\alpha_r(\sigma)T_{\sigma} \cdot (1+\bm{k}T_k)
=\sum_{w \preceq r}\alpha_r(w)T_{w} + \sum_{u \preceq
r}\alpha_r(u)T_{u}\cdot \bm{k}T_k.
\end{equation}
We will prove the desired result by applying \eqref{Tprod}, and
comparing the coefficients of $T_w$ on both sides of
\eqref{indLemma}.

\be
\item If $w, w s_k\npreceq r$ then we have
$\alpha_{r'}(w)=0=\alpha_{r}(w s_k)\cdot \bm{k}$. If $w\npreceq r$
and $w s_k \preceq r$ then $T_w$ can only be obtained by
   $T_{w s_k} \cdot \bm{k}T_k$, so we have
   $\alpha_{r'}(w)=\alpha_{r}(w s_k)\cdot \bm{k}$.

\item If $w \preceq r$ and $w s_k \npreceq r$, then there is no
   $u \preceq r$ such that $u s_k=w$. Hence $T_w$ can only be
   obtained by $T_{w} \cdot 1$, so we have
   $\alpha_{r'}(w)=\alpha_{r}(w)$.

\item If $w,w s_k \preceq r$ and $\ell(w s_k)=\ell(w)+1$, then
   $T_w \cdot \bm{k}T_k=\bm{k}T_{w s_k}$, and there is
   $u=w\cdot s_k \preceq r$ such that $T_{u} \cdot
   \bm{k}T_k=\bm{k}((q-1)T_{u}+q T_{u \cdot s_k})
   =\bm{k}((q-1)T_{w s_k}+q T_{w})$. Therefore we have
   $\alpha_{r'}(w)=
   \alpha_{r}(w)+ \alpha_{r}(u)\cdot \bm{k}q
   =\alpha_{r}(w)+ \alpha_{r}(ws_k)\cdot \bm{k}q$.

\item If $w,w s_k \preceq r$ and $\ell(w s_k)=\ell(w)-1$, then
   $T_w \cdot \bm{k}T_k=\bm{k}((q-1)T_{w}+qT_{w s_k})$, and there
   is $u=w\cdot s_k \preceq r$ such that $T_{u} \cdot
   \bm{k}T_k=\bm{k}T_{u \cdot s_k}=\bm{k}T_{w}$.  Therefore we
   have
\begin{eqnarray*}
\alpha_{r'}(w)&=& \alpha_{r}(w)+ \alpha_{r}(u)\cdot
\bm{k}+\alpha_{r}(w)\cdot \bm{k}(q-1)\\
&=&\alpha_{r}(w)\cdot q^k+ \alpha_{r}(ws_k)\cdot \bm{k}.
\end{eqnarray*}
\ee\qed

We also want to list the following result related to $\inv(w)$ and
$\inv(w s_k)$, which is frequently used in the proof of Theorem
\ref{mainthm}. The proof of this result is quite straightforward and
is omitted here.

\begin{lemma}\label{prop1}
For any permutation $w\in \sn$ and adjacent transposition $s_k,
1\leq k\leq n-1$, we have the following properties of the statistic
inv$_w$. \be
\item If $\ell(w s_k)=\ell(w)-1$, then
$$\inv_w(k)>\inv_w(k+1),\ \ \inv_{ws_k}(k)=\inv_w(k+1),$$
$$ \qquad \mbox{and }\ \inv_{ws_k}(k+1)=\inv_w(k)-1.$$

\item If $\ell(w s_k)=\ell(w)+1$ then
$$\inv_w(k)\leq \inv_w(k+1), \ \ \inv_{ws_k}(k)=\inv_w(k+1)+1,$$
$$ \qquad \mbox{and }\ \inv_{ws_k}(k+1)=\inv_w(k).$$
\ee
\end{lemma}

We need one more lemma before we prove the main theorem.

\begin{lemma}\label{inv}
Let $r$ be a tight sequence, and $w$ be a permutation such that $w
\preceq r$, then $\inv_w(i) \leq a_r(i)$ for any $i\geq 1$.
\end{lemma}
\proof We will prove this result by induction on the length of $r$.
The case for $r=(1)$ is trivial. Suppose the result is true for some
tight sequence $r$, and $r'=(r,k)$ is also a tight sequence. For any
permutation $w\preceq r'$, if $w\preceq r$, we have
\[\inv_w(i) \leq a_r(i) \leq a_{r'}(i),\ \ \forall i\geq 1. \]

If $w\npreceq r$, we will first prove that in this case
$a_r(k)=a_r(k+1)$. If $a_r(k)=0$, it is obvious that
$a_r(k)=a_r(k+1)=0$. If $a_r(k)>0$, then since $w\preceq r'$ but
$w\npreceq r$, there must be a $k+1$ to the right of the rightmost
$k$ in $r$. Suppose this $k+1$ occurs at the $j$th position in $r$.
We have that $(r_1,r_2,\ldots, r_{j-1},r_j)=(r_1,r_2,\ldots,
r_{j-1},k+1)$ is a tight 
sequence, with some permutation $w' \preceq (r_1,r_2,\ldots,
r_{j-1},r_j)$ but $w' \npreceq (r_1,r_2,\ldots, r_{j-1})$. Thus we
have $a_{(r_1,\ldots, r_{j-1})}(k+1)=a_{(r_1,\ldots,
r_{j-1})}(k)-1$. Therefore $a_r(k)=a_{(r_1,\ldots,
r_{j-1})}(k)=a_{(r_1,\ldots, r_{j-1})}(k+1)+1=a_r(k+1)$.

Moreover, since $w\npreceq r$, there exists a permutation $u
\npreceq r$ such that $w=us_k$ and $l(w)=l(u)+1$. It is
obvious that for $i \neq k,k+1$, we have $\inv_w(i)=\inv_u(i) \leq
a_r(i)=a_{r'}(i)$. Moreover, from Lemma~\ref{prop1} we know that
\[\inv_w(k)=\inv_u(k+1)+1 \leq a_r(k+1)+1=a_r(k)+1=a_{r'}(k),\]
and \[\inv_w(k+1)=\inv_u(k) \leq a_r(k)= a_r(k+1)=a_{r'}(k+1).\]

Hence the proof is complete. \qed

Now we are ready to prove the main theorem.

\noindent \emph{Proof of Theorem \ref{mainthm}.} The proof is by
induction on $l$, the length of the sequence $r$. It is trivial to
check that \eqref{main0} holds for $r=(1)$. Suppose that
\eqref{main0} holds for some tight sequence $r$, and $r'=(r, k)$ is
also a tight sequence.  We want to prove that \eqref{main0} also
holds for $r'$. The case when $k
> \max\{r\}$ is trivial, so from now on we will assume that
$\max\{r\}=\max\{r'\}=n-1$.

For any $i\ne k,k+1$ ($1\leq i \leq n-1$), we have
$a_{r'}(i)=a_{r}(i)$, and $\inv_w(i)=\inv_{ws_k}(i)$. Therefore
\beas \max \{a_{r'}(i-1)-1,\inv_w(i)\} & = & \max
\{a_{r}(i-1)-1,\inv_w(i)\}\\
& = & \max \{a_{r}(i-1)-1,\inv_{w s_k}(i)\} \eeas holds for any
$i\ne k,k+1$ ($2\leq i \leq n-1$). Hence we only need to concentrate
on the values of $\max \{a_{r'}(i-1)-1,\inv_w(i)\}$ for $i=k,k+1$.
(When $k=1$, we only consider $\max \{a_{r'}(1)-1,\inv_w(2)\}$.)

Next we will prove that
$\alpha_{r'}(w)=\prod_{i=2}^{n-1}\bm{i}^{\max\{a_{r'}(i-1)-1,
   \inv_{w}(i)\}}$ for any $w \npreceq r'$ according to the four cases in Lemma \ref{lem0},
and we will frequently use Lemma~\ref{prop1}.

\be
\item Let $w\npreceq r$. In this case
   $\inv_w(k)=a_{r}(k)+1=a_{r'}(k)\leq n-k$.  Since $r, r'$ are both
   tight sequences we have $a_{r'}(k-1)=a_{r}(k-1)=a_{r}(k)+1$.
   Moreover, since $\inv_{ws_k}(k)=\inv_w(k+1)<\inv_w(k)=a_{r}(k)+1$,
   we have
$$\max \{a_{r'}(k-1)-1,\inv_{w}(k)\}=a_{r}(k)+1=\max \{a_{r}(k-1)-1,\inv_{ws_k}(k)\}+1.$$

Since $\inv_{ws_k}(k+1)=\inv_w(k)-1=a_{r}(k)$ and
$a_{r'}(k)=a_{r}(k)+1$, we have
$$\max \{a_{r'}(k)-1,\inv_{w}(k+1)\}=\max
 \{a_{r}(k)-1,\inv_{ws_k}(k+1)\}=a_{r}(k).$$

Hence we conclude that
  \beas \alpha_{r'}(w)=\alpha_{r}(ws_k)\cdot \bm{k} & = &
\prod_{i=2}^{n-1}\bm{i}^{\max \{a_{r}(i-1)-1,\inv_{ws_k}(i)\}}\cdot
\bm{k}\\
& = & \prod_{i=2}^{n-1}\bm{i}^{\max \{a_{r'}(i-1)-1,\inv_{w}(i)\}}.
\eeas

\item Let $w\preceq r$ and $w s_k \npreceq r$. In this case we have
   $\ell(w s_k)=\ell(w)+1$ and $\inv_w(k)= a_{r}(k)$.

Since $\inv_w(k+1)\geq \inv_w(k)=a_{r}(k)$ and
$a_{r'}(k)=a_{r}(k)+1$, we have
$$\max \{a_{r}(k)-1,\inv_w(k+1)\}=\max \{a_{r'}(k)-1,\inv_w(k+1)\}.$$
It follows that $\alpha_{r'}(w)=\alpha_{r}(w)
=\prod_{i=2}^{n-1}\bm{i}^{\max \{a_{r'}(i-1)-1,\inv_{w}(i)\}}$.

\item Let $w,w s_k \preceq r$ and $\ell(w s_k)=\ell(w)+1$.
Since $\inv_w(k)<a_{r}(k)$, $\inv_{ws_k}(k)\leq a_{r}(k)$ and
$a_{r}(k-1)-1\geq a_{r}(k)$, we have
$$\max \{a_{r}(k-1)-1,\inv_w(k)\}=\max \{a_{r}(k-1)-1,\inv_{ws_k}(k)\}=a_{r}(k-1)-1.$$

Since $\inv_w(k+1)=\inv_{ws_k}(k)-1 \leq a_{r}(k)-1$, and $\inv_{w
   s_k}(k+1)=\inv_w(k)<a_{r}(k)$, we have
$$\max \{a_{r}(k)-1,\inv_{w}(k+1)\}=\max \{a_{r}(k)-1,\inv_{ws_k}(k+1)\}=a_{r}(k)-1.$$
Hence $\alpha_r(w)=\alpha_r(ws_k)$. Therefore we have \beas
   \alpha_{r'}(w) = \alpha_{r}(w)+ \alpha_{r}(ws_k)\cdot \bm{k}q
    & = & \alpha_{r}(w)\bm{(k+1)}\\ & = &
   \prod_{i=2}^{n-1}\bm{i}^{\max \{a_{r}(i-1)-1,\inv_{w}(i)\}}\cdot
   \bm{(k+1)}.
\eeas Moreover, since $\max \{a_{r'}(k)-1,\inv_{w}(k+1)\}=a_{r}(k)=
\max \{a_{r}(k)-1,\inv_{w}(k+1)\}+1$, we have
$\alpha_{r'}(w)=\prod_{i=2}^{n-1}\bm{i}^{\max\{a_{r'}(i-1)-1,\inv_{w}(i)\}}$.

\item Let $w,w s_k \preceq r$ and $\ell(w s_k)=\ell(w)-1$. In this
   case $\inv_w(k)\leq a_{r}(k)$.

Since $\inv_{ws_k}(k)=\inv_w(k+1)<\inv_w(k)\leq a_{r}(k)$ and
$a_{r}(k-1)-1\geq a_{r}(k)$, we have
$$\max \{a_{r}(k-1)-1,\inv_{w}(k)\}=\max
\{a_{r}(k-1)-1,\inv_{ws_k}(k)\}=a_{r}(k-1)-1.$$ Since
$\inv_{ws_k}(k+1)=\inv_w(k)-1$, we have
$$\max \{a_{r}(k)-1,\inv_{w}(k+1)\}=\max
\{a_{r}(k)-1,\inv_{ws_k}(k+1)\}=a_{r}(k)-1.$$ Hence
$\alpha_r(w)=\alpha_r(ws_k)$. Therefore we have \beas \alpha_{r'}(w)
=\alpha_{r}(w)\cdot q^k+ \alpha_{r}(ws_k)\cdot \bm{k} & = &
\alpha_{r}(w)\cdot \bm{(k+1)}\\ & = & \prod_{i=2}^{n-1}\bm{i}^{\max
\{a_{r}(i-1)-1,\inv_{w}(i)\}}\cdot \bm{(k+1)}. \eeas Moreover, since
$\max \{a_{r'}(k)-1,\inv_{w}(k+1)\}=a_{r}(k) =\max
\{a_{r}(k)-1,\inv_{w}(k+1)\}+1$, we have
$\alpha_{r'}(w)=\prod_{i=2}^{n-1}\bm{i}^{\max\{a_{r'}(i-1)-1,\inv_{w}(i)\}}$.
\ee

Hence the proof is complete. \qed

We can use Theorem~\ref{mainthm} and its proof to compute
$\alpha_r(w)$ for certain sequences $r$ that are not tight
sequences.

\begin{corollary} \label{coro1}
Let $r$ be a sequence of positive integers, and $\max\{r\}=n-1$. If
$r$ has the prefix $\rho_n=(1,2,1,3,2,1\ldots, n,n-1,\ldots,1)$,
then we have
\begin{equation} \label{eq:alrw}
\alpha_r(w)=\prod_{i=2}^{n-1}
\bm{i}^{\max\{a_{r}(i-1)-1,\inv_{w}(i)\}}.
\end{equation}
\end{corollary}

\proof We will prove equation~\eqref{eq:alrw} by induction on the
length of $r$. Since $\rho_n$ is a tight sequence, from Theorem
\ref{mainthm} we know that the result holds for $r=\rho_n$. Next
assume the result for $r$ and let $r'=(r,k)$ with $1 \leq k \leq
n-1$. We do an induction similar to what we did in the proof of
Theorem \ref{mainthm}.  Since $r$ has the prefix $\rho_n$, it
follows that for any $w \in \sn$, $w, w s_k \preceq r$. Therefore
only cases 3 and 4 will occur. Moreover, since $a_r(k-1)\geq
n-(k-1)$, $a_r(k)\geq n-k$ and $a_{r'}(k)=a_{r}(k)+1$, we have
$$\max \{a_{r}(k-1)-1,\inv_{w}(k)\}=\max
\{a_{r}(k-1)-1,\inv_{ws_k}(k)\}=a_{r}(k-1)-1,$$
$$\max \{a_{r}(k)-1,\inv_{w}(k+1)\}=\max
\{a_{r}(k)-1,\inv_{ws_k}(k+1)\}=a_{r}(k)-1,$$ and
$$\max \{a_{r'}(k)-1,\inv_{w}(k+1)\}=a_{r}(k)
=\max \{a_{r}(k)-1,\inv_{w}(k+1)\}+1.$$ Hence for both case 3 and 4
we have
$\alpha_{r'}(w)=\prod_{i=2}^{n-1}\bm{i}^{\max\{a_{r'}(i-1)-1,\inv_{w}(i)\}}$.
\qed.

Note that $r$ is a reduced decomposition of $w \in \sn$ if and only
if the reverse of $r$ is a reduced decomposition of $w^{-1}$. Thus
we have the following result.
\begin{corollary} \label{coro2}
Let $r$ be a sequence of positive integers, and $\max\{r\}=n-1$. If
\be
\item $r$ is the reverse of a tight sequence, or
\item $r$ has suffix
$\rho_n=(1,2,1,3,2,1\ldots, n,n-1,\ldots,1)$, \ee then for any $w
\in \sn$ and $w \preceq r$, we have
\begin{equation}
\alpha_r(w)=\prod_{i=2}^{n-1}
\bm{i}^{\max\{a_{r}(i-1)-1,\inv_{w^{-1}}(i)\}}.
\end{equation}
\end{corollary}

\textsc{Note.} If a sequence $r'$ is obtained from $r$ by
transposing two adjacent terms that differ by at least 2, then
$Q(r)=Q(r')$, so $\alpha_w(r)=\alpha_w(r')$. Thus our results extend
to sequences that can be obtained from those of
Theorem~\ref{mainthm}, Corollary~\ref{coro1}, and
Corollary~\ref{coro2} by applying such ``commuting transpositions''
to $r$.

\section
{A Connection with Random Walks on $\sn$} There is a huge literature
on random walks on $\sn$, e.g., \cite{PDiaconis}. Our results can be
interpreted in this context. First consider the case $q=1$. In this
case the Hecke algebra $\hnq$ reduces to the group algebra $\rr\sn$
of $\sn$, and the generator $T_i$ becomes the adjacent transposition
$s_i$. Thus
   $$ Q(r)_{q=1}=(1+r_1s_{r_1})(1+r_2s_{r_2})\cdots(1+r_ls_{r_l}). $$
We normalize this expression by dividing each factor $1+r_is_i$ by
$1+r_i$. Write
    $$ D_j=(1+js_j)/(1+j), $$
and set
   $$ \widetilde{Q}(r) = D_{r_1}D_{r_2}\cdots D_{r_l}. $$
If $P$ is a probability distribution on $\sn$, then let $\sigma_P =
\sum_{w\in\sn} P(w)w\in\rr\sn$. If $P'$ is another probability
distribution on $\sn$, then $\sigma_P\sigma_{P'} =\sigma_{P*P'}$ for
some probability distribution $P*P'$, the \emph{convolution} of $P$
and $P'$. It follows that $\widetilde{Q}(r) = \sigma_{P_r}$ for some
probability distribution $P_r$ on $\sn$. Theorem~\ref{mainthm} gives
(after setting $q=1$ and normalizing) an explicit formula for the
distribution $P_r$, i.e., the values $P_r(w)$ for all $w\in \sn$.
Note in particular that if $r$ is the standard tight sequence
$\rho_n = (1,2,1,3,2,1,4,3,2,1,\dots,n-1,n-2,\dots,1)$, then from
equation~\eqref{eq:qunif} we get
   $$ \widetilde{Q}(\rho_n) = \frac{1}{n!}\sum_{w\in\sn} w = \sigma_U,
   $$
where $U$ is the uniform distribution on $\sn$.  (We have been
informed by Alexander Molev that an equivalent result was given by
Jucys \cite{jucys} in 1966. We have also been informed by Persi
Diaconis that this result, and similar results for some other
groups, were known by him and Colin Mallows twenty years ago.)  It
is not hard to see directly why we obtain the uniform distribution.
Namely, start with any permutation $w=w_1\cdots w_n\in\sn$. Do
nothing with probability 1/2 or apply $s_1$ (i.e., interchange $w_1$
and $w_2$) with probability 1/2, obtaining $y_1y_2 w_3\cdots w_n$.
Thus $y_2$ is equally likely to be $w_1$ or $w_2$. Now either do
nothing with probability 1/3 or apply $s_2$ with probability 2/3,
obtaining $y_1 z_2 z_3 w_4 \cdots w_n$. Then $z_3$ is equally likely
to be $w_1, w_2$ or $w_3$.  Continue in this way, applying
$s_3$,\dots, $s_{n-1}$ at each step or doing nothing, with
probability $1/(i+1)$ of doing nothing at the $i$th step, obtaining
$d_1\cdots d_n$. Then $d_n$ is equally likely to be any of
$1,2,\dots,n$. Now apply $s_1, s_2, \dots, s_{n-2}$ or do nothing as
before, obtaining $e_1 \cdots e_n$. The last element $e_n$ has never
switched, so $e_n=d_n$, and now $e_{n-1}$ is equally likely to be
any element of $\{1,2,\dots,n\} -\{d_n\}$. Continue as before with
$s_1,\dots,s_{n-3}$, then $s_1,\dots,s_{n-4}$, etc., ending in $s_1,
s_2, s_1$, at which point we obtain a uniformly distributed random
permutation.

Now consider the situation for $\hnq$. If $P$ is a probability
distribution on $\sn$ then write $\tau_P=\sum_{w\in\sn} P(w)T_w
\in\hnq$. If $P'$ is another probability distribution on $\sn$, then
in general it is not true that $\tau_P\tau_P'=\tau_R$ for some
probability distribution $R$. A probabilistic interpretation of
Theorem~\ref{mainthm} requires the use of a Markov chain. Let $0< q<
1$. Note that from equation~\eqref{Tprod} we have
   $$ T_w(1+\bm{k}T_k) = \left\{ \begin{array}{rl}
       T_w + \bm{k}T_{ws_k}, & \ell(ws_k)=\ell(w)+1\\
          q^k T_w+q\bm{k}T_{ws_k}, & \ell(ws_k)=\ell(w)-1.
         \end{array} \right. $$
Divide each side by $1+\bm{k}$. Let $w=w_1 \cdots w_n$. We can then
interpret multiplication of $T_w$ by $(1+\bm{k}T_w)/(1+\bm{k})$ as
follows. If $w_k<w_{k+1}$ then transpose $w_k$ and $w_{k+1}$ with
probability $\bm{k}/(1+\bm{k})$, or do nothing with probability
$1/(1+\bm{k})$. If $w_k>w_{k+1}$, then transpose $w_k$ and $w_{k+1}$
with probability $q\bm{k}/(1+\bm{k})$, or do nothing with
probability $q^k/ (1+ \bm{k})$. Since
   $$ \frac{q\bm{k}}{1+\bm{k}}+\frac{q^k}{1+\bm{k}}<1, $$
we have a ``leftover'' probability of
$(1-(q\bm{k}+q^k)/(1+\bm{k}))$. In this case the process has failed
and we should start it all over. Let us call this procedure a
\emph{$k$-step}.

If $r=(r_1,\dots,r_l)$ is a tight sequence, then begin with the
identity permutation and apply an $r_1$-step, $r_2$-step, etc. If we
need to start over, then we again begin with the identity
permutation and apply an $r_1$-step, $r_2$-step, etc. Eventually
(with probability 1) we will apply $r_i$-steps for all $1\leq i\leq
l$, ending with a random permutation $v$. In this case,
Theorem~\ref{mainthm} tells us the distribution of $v$, namely, the
probability of $v$ is
   $$ P(v) = \frac{\alpha_r(v)}{\prod_{i=1}^l (1+\bm{r_i})}. $$
In particular, if $r$ is the standard tight sequence $\rho_n$, then
$v$ is uniformly distributed.

\begin{example}
Start with the permutation 123 and $r=\rho_3=(1,2,1)$. Let us
calculate by ``brute force'' the probability $P=P(123)$ that
$v=123$. There are three ways to achieve $v=123$.
   \be \item[(a)] Apply a 1-step, a 2-step, and a 1-step, doing nothing
   each time. This has probability $(1/2)(1/(2+q))(1/2) = 1/4(2+q)$.
    \item[(b)] Apply a 1-step and switch. Apply a 2-step and do
      nothing. Apply a 1-step and switch. This has probability
      $q/4(2+q)$.
    \item[(c)] Apply a 1-step and switch. Apply a 2-step and do
      nothing. Try to apply a 1-step but go back to the beginning,
      after which we continue the process until ending up with
      123. This has probability
        $$ \frac 12\frac{1}{2+q}(1-q)P = \frac{P(1-q)}{2(2+q)}. $$
Hence
   $$ P =\frac{1}{4(2+q)}+\frac{q}{4(2+q)} +\frac{P(1-q)}
          {2(2+q)}. $$
Solving for $P$ gives (somewhat miraculously!) $P=1/6$. Similarly
for all other $w\in\mathfrak{S}_3$ we get $P(w)=1/6$.
   \ee
\end{example}

\textsc{Note.} A probabilistic interpretation of certain Hecke
algebra products different from ours appears in a paper by Diaconis
and Ram \cite{d-r}.

\end{document}